\documentclass[12pt,leqno]{article} 
\usepackage{amsmath,amsfonts,amscd,Martintxo}
\usepackage[colorlinks]{hyperref}

\title{\Large\bf  {De Rham intersection cohomology for 
general perversities.}}

\author{
Martintxo Saralegi-Aranguren\thanks{EA 2462 Laboratoire de 
Mathématiques de Lens. 
Faculté Jean Perrin.
Université d'Artois.   Rue Jean Souvraz SP 18.   62 307 Lens Cedex - 
France.   
{\sl saralegi@euler.univ-artois.fr}. }
\\ {\small Université d'Artois }
}
\date{}
\begin{document}  
\maketitle

\thispagestyle{empty}

\begin{abstract} 
\nt For a stratified pseudomanifold $X$,  we have the de Rham Theorem 
$
\lau{\IH}{*}{\per{p}}{X} = 
\lau{\IH}{\per{t} - \per{p}}{*}{X},
$
for a perversity $\per{p}$ verifying $\per{0} \leq \per{p} \leq 
\per{t}$, where $\per{t}$ denotes the top perversity.
We extend  this result to any perversity $\per{p}$. In the direction 
cohomology $\mapsto$ homology, we obtain the 
isomorphism
$$
\lau{\IH}{*}{\per{p}}{X} = 
\lau{\IH}{\per{t} -\per{p}}{*}{X,\ib{X}{\per{p}}},
$$
where 
$
{\displaystyle \ib{X}{\per{p}} = \bigcup_{ S \preceq S_{1} \atop
\per{p} (S_{1})< 0}S = \bigcup_{ \per{p} (S)< 0}\overline{S}.}
$
In the direction 
homology $\mapsto$ cohomology, we obtain the 
isomorphism
$$
\lau{\IH}{\per{p}}{*}{X}=\lau{\IH}{*}{\max ( \per{0},\per{t} -\per{p})}{X}.
$$
In our paper stratified pseudomanifolds with 
one-codimensional strata are allowed.
 \end{abstract}

 Roughly speaking, a stratified pseudomanifold $X$ is a family ${\cal S}_X$ of smooth 
 manifolds (strata) assembled in a conical way.
 A (general) perversity $\per{p}$ associates an integer to each of the strata of 
 $X$ (see \cite{M}).
 The classical perversities (see \cite{GM1}, \cite{K}, \cite{F}, 
 \ldots ) are filtration-preserving, that is, they 
 verify:
 $$
 S_{1},S_{2} \in {\cal S}_X \hbox{ with } \dim S_{1} = \dim S_{2} 
 \Rightarrow \per{p}(S_{1}) = \per{p}(S_{2}).
 $$
 The zero-perversity, defined by $\per{0}(S) = 0$, and the top 
 perversity, defined by $\per{t}(S) = \codim_{X} S -2$, are classical 
 perversities.

 The singular intersection homology 
 $ \lau{\IH}{\per{p}}{*}{X} $ was introduced by 
 Goresky-MacPherson in \cite{GM2} (see also \cite{K}). It is a 
 topological invariant of the stratified pseudomanifold when the 
 perversity satisfies some monotonicity conditions (see \cite{GM1}, 
 \cite{K}, \ldots ). In particular, we need $\per{0} \leq \per{t}$ and 
 therefore $X$
  does not possess any one-codimensional strata.
  Recently, a more general result has been obtained in \cite{F} where 
  one-codimensional strata are allowed. In all these cases, the 
  perversities are classical.
  
  The de Rham intersection cohomology   $\lau{\IH}{*}{\per{p}}{X}$ was 
  also introduced by 
  Goresky-MacPherson (see \cite{Bry}). It requires the existence of a 
  Thom-Mather neighborhood system. Other versions exist, but always an extra 
  datum is needed in order to define this cohomology: a   
  Thom-Mather neighborhood system
  (\cite{Bry}, \cite{BL} \ldots , ) a riemannian 
  metric (\cite{ChGM},  \cite{N}, \cite{BHS2}, \ldots), a PL-structure (\cite{BGM}, 
  \cite{CHS}, \ldots ), a blow up (\cite{BHS}, \ldots), etc.

 The perverse de Rham Theorem 
\begin{equation}
    \label{DR}
 \lau{\IH}{*}{\per{p}}{X} =\lau{\IH}{\per{t} - \per{p}}{*}{X},
\end{equation}
 relates the intersection 
 homology  with the intersection cohomology. 
First  it was proved  by  Brylinski in \cite{Bry} and after in the above references. 
The involved perversities are classical perversities verifying some 
monotonicity conditions. Moreover, the perversity $\per{p}$ must lie 
between $\per{0}$ and $\per{t}$, which exclude the existence of 
one-codimensional strata on $X$.

The first proof of the de Rham Theorem for the general perversities 
has been
given by the author in \cite{S}  using the integration ${\displaystyle 
\int}$ of 
 differential forms on simplices. Unfortunately, there is a mistake in 
 the  statement of Proposition 2.1.4 and Proposition 2.2.5: 
     the 
 hypothesis $\per{p} \leq \per{t}$ must be added.
 As a consequence, the main result\footnote{The statement 
 of the second main result (the Poincaré 
 Duality: Theorem 4.2.7) does not need any modification since 
 Propositions 2.1.4 and 2.2.5 are not used for its proof.} of \cite{S} 
 (de Rham Theorem 4.1.5) is valid for 
 a general perversity $\per{p}$ \underline{verifying the condition $\per{0} \leq \per{p} \leq 
 \per{t}$}. In particular, we have \refp{DR} for a general perversity $\per{p}$ 
 with $\per{0} \leq \per{p} \leq 
 \per{t}$. Notice that the one-codimensional strata are not allowed.

\bigskip

 In this work we prove a de Rham Theorem for  any 
general perversity $\per{p}$\footnote{The one codimensional strata are 
finally allowed! }. The formula \refp{DR} changes! We obtain that, in the direction cohomology $\mapsto$ homology, 
 the integration ${\displaystyle \int}$ induces the isomorphism
 $
 \lau{\IH}{*}{\per{p}}{X} = 
\lau{\IH}{\per{t} -\per{p}}{*}{X,\ib{X}{\per{p}}},
 $
 where 
 $
 {\displaystyle \ib{X}{\per{p}} = \bigcup_{ S \preceq S_{1} \atop
 \per{p} (S_{1})< 0}S = \bigcup_{ \per{p} (S)< 0}\overline{S}.}
 $
 (cf. Theorem \ref{Rham}).
 In the direction 
 homology $\mapsto$ cohomology, we have the 
 isomorphism
 $
 \lau{\IH}{\per{p}}{*}{X}
 =
\lau{\IH}{*}{\max ( \per{0},\per{t} -\per{p})}{X},
 $
 (cf. Corollary \ref{Rhambis}).

 We end the work by noticing that the Poincaré Duality of \cite{Bry} and 
 \cite{S} is still valid in our context.

We sincerely thank the anonymous referee for his/her work on this 
paper and the helpful comments provided in the review. They have 
definitely helped  to improve the paper, to make it clearer in the 
objectives and easier to read.
 
 \section{Stratified spaces and unfoldings.}

 We present the geometrical framework of this work, that is, the stratified 
 pseudomanifolds and the unfoldings. For a more complete study of these notions, we refer the reader to, for 
 example,    \cite{GM2} and \cite{S}.

 In the sequel, any manifold is connected, second 
 countable, Haussdorff,   
 without boundary and smooth 
 (of class $C^\infty$).

 \prg {\bf Stratifications.} A {\em stratification} of
 a  paracompact  space $X$ is a locally
 finite partition ${\cal S}_X$ of $X$ into disjoint smooth manifolds,
 called {\em strata}, such that

 $$
 S\cap \overline{S'}\not= \emptyset 
 \Longleftrightarrow S \subset \overline{S'} ,
 $$
 written $S \preceq S'$.
 Notice that $({\cal S}_X, \preceq)$ is a partially ordered set. 
 
 We  say that $X$ is a  {\em stratified space}.
 The {\em depth} of $X$, denoted $\depth X$,  is the length 
 of the maximal chain contained in $X$.  It is always finite because of the locally 
finiteness of ${\cal S}_X$. The minimal (resp. maximal) strata are the closed (resp. open)
strata.
 The open strata are the {\em regular strata} and the other ones  are the 
 {\em singular strata}. We shall denote ${\cal S}^{^{sing}}_X$ the family of 
 singular strata. The union $\Sigma_X$ of singular strata is the {\em 
 singular part}, 
 which is a closed subset. The {\em regular part} $X - \Sigma_X$ is an open 
 dense subset.  We require the regular strata to have the same dimension, 
denoted $\dim X$.

 For each $i \in \{ -1, 0 , \ldots, \dim X\}$ we  consider the saturated 
 subset:
 $$
 \ib{X}{i} = \bigcup \{ S \in {\cal S}_X \ / \ \dim S \leq i \}. 
 $$
 This gives the {\em filtration} $\cal{F}_{X}$:
\begin{equation}
    \label{FF}
 \ib{X}{\dim X}\supset \ib{X}{\dim X-1} \supset\ib{X}{1}\supset \ib{X}{0} \supset 
 \ib{X}{-1} =\emptyset.
\end{equation}

 The main example of a stratified space is given by the following conical 
 construction.
 Consider  a compact stratified space $L$ and let $cL$ be the {\em cone} of 
 $L$, that is  $cL = L \times [0,1[\Big/ L \times \{ 0 \}$. 
 The points of $cL$ are denoted by $[x,t]$.
 The {\em vertex} of the cone is the point  $\vartheta = [x,0]$.
 This cone is 
 naturally endowed with the following stratification:
 $$
 {\cal S}_{cL} =
 \{ \{\vartheta \} \} \cup \left\{ S\times ]0,1[ \ / \ S \in {\cal S}_L\right\}.
 $$
 For the filtration $\mathcal{F}_{cL}$ we have:
 $$
 \left( cL\right)_i =
 \left\{ 
 \begin{array}{ll}
 \{ \vartheta \} & \hbox{if $i=0$}\\
 cL_{i-1} & \hbox{if $i>0$}
 \end{array}.
 \right.
 $$
 Notice that $\depth cL = \depth L +1$.

 The {\em 
 canonical stratification} of a manifold $X$
 is the family
 ${\cal S}_X $ formed by the connected components of $X$. The 
 filtration contains just one non-empty element: $X_{\dim X}$.

 A continuous map (resp. homeomorphism) $f \colon Y \to
 X$ between two stratified spaces is a stratified morphism (resp.
 {\em isomorphism}) if it sends  the strata of  $Y$  to the strata of
 $X$ smoothly (resp.
 diffeomorphically). 

 \prg {\bf Stratified pseudomanifolds}. A stratified space $X$ is a {\em 
 stratified pseudomanifold} when it possesses a conical local structure. 
 More explicitly, when for each point $x$ of a singular  stratum $S$ of $X$ there exists a
 stratified isomorphism
 $
 \phii \colon U \TO \R^n \times cL_S,
 $
 where

 \Zati $U \subset X$ is an open neighborhood of $x$ endowed with the 
 induced stratification,

 \zati $L_S$  is a compact stratified space, called {\em link} of $S$,

 \zati $\R^n \times cL_S$ is endowed with the stratification 
 $
 \left\{ \R^n \times \{ \vartheta \} \right\} \cup \left\{ \R^n \times 
 S' \times ]0,1[\ / \ S' \in {\cal S}_{L_S} \right\}$, and

 \zati $\phii(x) = (0,\vartheta)$.

\smallskip

 The couple $(U,\phii)$ is a {\em chart} of $X$ containing $x$. An {\em atlas} ${\cal A}$ is 
 a family of charts covering $X$. A stratified 
 pseudomanifold is {\em normal} when all the links are connected. 
 Notice that in this case each link is a connected normal stratified 
 pseudomanifold.

 \prg{\bf Unfoldings.}
 Consider  a stratified pseudomanifold $X$. A continuous
 map  $ {\cal L}
 \colon \wt{X} \to X$, where  $\wt{X}$ is a (not
 necessarily connected) manifold,  is an {\em unfolding}
 if  the two following conditions hold: 
 \begin{itemize}

  \item[1.] The restriction  $ {\cal L}_{X} \colon  {\cal
 L}_{X}^{-1} (X - \Sigma_X ) 
 \longrightarrow X - \Sigma_X $ is a local diffeomorphism.

 \item[2.]   There exist a family of unfoldings  $ \{
 {\cal L}_{L_S}
 \colon  \wt{L_S} \to L_S\}_{S\in {\cal S}_X^{^{sing}}} $
 and an atlas ${\cal A}$ of
 $X$ such that for  each chart $(U,\phii)\in {\cal A}$ 
 there exists a commutative diagram
$$
 \begin{CD}
     \R^{n} \times \wt{L_S}
     \times ]-1,1[@>\wt{\phii}>>{\cal L}_X^{-1}(U)\\
     @VQVV @V\mathcal{L}_{X}VV\\
     \R^{n}  \times c L_S@>\phii>> U\\
     \end{CD}
     $$
 where

 \Zati   $\wt{\phii}$ is a diffeomorphism  and  

 \zati $Q(x,\wt{\zeta},t) = (x, \left[ {\cal L}_{L_S}(\wt{\zeta}),|t|
 \right])
 $.
\end{itemize}

  We  say that $X$ is an {\em unfoldable pseudomanifold}. 
 This definition makes sense because it is made by
 induction on $\depth X$.   When $\depth X =0$ then
 ${\cal L}_X$ is just a local diffeomorphism. 
 For any
 singular stratum
 $S$ the
 restriction ${\cal L}_X \colon  {\cal L}_X^{-1}(S)\to
 S$ is a fibration with  fiber $\wt{L}_S$.     
 The {\em canonical unfolding } of the cone $cL_S$ is  the map ${\cal L}_{cL_S} \colon 
 \wt{cL_S}  = \wt{L_S}   \times ]-1,1[ \to cL_S$ 
 defined by ${\cal L}_{cL_S}(\wt{\zeta},t) = \left[ {\cal L}_{L_S}(\wt{\zeta}),|t|
 \right].
 $

 From now on, $(X,\cal{S}_{X})$ is a stratified pseudomanifold
 endowed with an unfolding  $ {\cal L}_X \colon \wt{X}
 \to X$.
 
 \prg {\bf Bredon's Trick.}
The typical result we prove in  this work looks like the following affirmation:
\begin{center}
``The differential 
operator $f \colon \coho{A}{*}{X} \to \coho{B}{*}{X}$, defined between two differential 
complexes on $X$, induces an isomorphism in cohomology.''
\end{center}
First we  prove  this assertion for charts. The passing from  
local to global  can be done  using different  tools.  
For 
example, 
\begin{itemize}
    \item[-] Axiomatic presentation of the intersection homology 
    (the most employed: \cite{GM1}, \cite{Bry}, \cite{BHS}, \ldots),
    \item[-] Uniqueness of the minimal stratification (used in  
    \cite{K}),
    \item[-] The generalized Mayer-Vietoris principle of
\cite[Chapter II]{BT} (used in  \cite{S}). 
    \item[-] The 
Bredon's trick of \cite[page 289]{Br}.
    \end{itemize}
In this work we choose the last one, maybe the less technical. 
The exact statement is the following:

\bL
    \label{Brt}
    Let $Y$ be a paracompact topological space
and  let $\{ U_\alpha\}$ be an open covering, closed for 
finite intersection. Suppose that $Q(U)$ is a statement about open subsets 
of $Y$, satisfying the following three properties:
\begin{itemize}
    \item[]
    \begin{itemize}
\item[(BT1)] $Q(U_\alpha)$ is true for each $\alpha$;
\item[(BT2)] $Q(U)$, $Q(V)$ and $Q(U\cap V )$ $\Longrightarrow$ $Q(U\cup V 
)$, where $U$ and $V$ are open subsets of $Y$; 
\item[(BT3)] $Q(U_i) \Longrightarrow Q\left({\displaystyle  \ \bigcup_i 
}U_i\right)$, 
where $\{ U_i\}$ is a disjoint family  of open subsets of $Y$.
\end{itemize}
\end{itemize}
Then $Q(Y)$ is true.
\eL
 \section{Intersection homology.}
 The intersection homology was introduced by Goresky-MacPherson in 
 \cite{GM1}, \cite{GM2}. Here we use  the singular intersection 
 homology of \cite{K}.

 \prg{\bf Perversity.}  Intersection cohomology requires the 
 definition of a  perversity parameter $\per{p}$. 
It associates an integer to each singular
 stratum of X, in other words, a  {\em perversity} is  a map $\per{p}\colon {\cal 
 S}^{^{sing}}_X \to \Z$. The {\em zero perversity}  
 $\per{0}$  is defined by $\per{0}(S) =0$. The {\em top perversity} $\per{t}$ is defined 
 by $\per{t}(S) = \codim_X S -2$. Notice that the condition $\per{0} 
 \leq \per{t}$ implies $\codim_X S \geq 2$ for each singular stratum 
 $S$, and therefore,   the one-codimensional strata are not allowed.

 The classical perversities (cf.  \cite{GM1}, \cite{GM2}, \ldots ), 
 the loose perversities (cf. \cite{K}), the superperversities (cf. 
 \cite{CS}, \cite{F}, \ldots ), \ldots  are  filtration-preserving map: $\per{p}(S) = 
 \per{p}(S')$ if $\dim S = \dim S'$.
 They also verify a monotonicity condition and, for some of them, the 
 one-codimensional strata are avoided. For such perversities the 
 associated intersection cohomology is a topological invariant.
 
 In our case, the perversities are stratum-preserving without any 
 constraint. Of course, 
 the topological invariance is lost. But we prove that we have a 
 de Rham duality (between the intersection homology and the intersection 
 cohomology) and the Poincaré duality.

 \medskip 

 We fix   a perversity $\per{p}$.  The homologies and 
 the cohomologies we 
 use in this work are with coefficients in $\R$.

 \prg{\bf Intersection homology. First approach.}  
 A singular simplex $\sigma \colon \Delta \to X$ is a $\per{p}$-{\em allowable}
 simplex  if

 \begin{itemize}
 \item[(All)] $\sigma^{-1}(S) \subset (\dim \Delta - 2 - \per{t}(S) + 
 \per{p}(S))$-skeleton of $\Delta$, for each singular stratum $S$.
 \end{itemize}
 \smallskip

 \nt A singular chain $\xi = \sumas{j=1}{m} r_{j} 
 \sigma_{j} \in \homo{S}{*}{X}$ is $\per{p}$-{\em allowable}  if each singular simplex $\sigma_{j}$ is 
 $\per{p}$-allowable. 
  The family of 
 $\per{p}$-allowable chains is a graded vector space denoted by 
 $\lau{AC}{\per{p}}{*}{X}$. 
 The associated differential complex is the complex of
 {\em $\per{p}$-intersection
 chains}, that is, 
 $
 \lau{SC}{\per{p}}{*}{X}= \lau{AC}{\per{p}}{*}{X}
 \cap \partial^{-1}\lau{AC}{\per{p}}{*-1}{X}.
 $
 Its homology $
 \lau{\IH}{\per{p}}{*}{X} = \homo{H}{*}{\lau{SC}{\per{p}}{\cdot}{X}}
 $
 is the $\per{p}$-{\em intersection homology} of $X$.
 This is the approach of \cite{GM2}. The intersection homology verifies two important 
 computational properties: Mayer-Vietoris and the product formula 
 $\lau{\IH}{\per{p}}{*}{\R \times X}  = \lau{\IH}{\per{p}}{*}{X} $ (see 
 \cite{GM1}, \cite{GM2}, \cite{S}, \ldots). 
 
 The usual local calculation is  the following (cf. \cite{GM2}, 
 see also \cite{K}). It corrects Proposition 2.1.4 of \cite{S}.
 \bP
 \label{cono1}
 Let $L$ be a compact stratified pseudomanifold. Then
 $$\lau{\IH}{\per{p}}{i}{cL} =
 \left\{\begin{array}{ll}
 \lau{\IH}{\per{p}}{i}{L} & \hbox{if $i \leq  \per{t}(\vartheta)- \per{p}(\vartheta)$} 
 \\[,2cm]
 0& \hbox{if $0 \ne  i \geq 1+ \per{t}(\vartheta)- \per{p}(\vartheta)$}
 \\[,2cm]
 \R & \hbox{if $0 = i \geq 1 + \per{t}(\vartheta)- \per{p}(\vartheta)$}.
 \end{array}
 \right.
 $$
 \eP
 \pro
 For $i \leq 1 +\per{t}(\vartheta) - \per{p}(\vartheta)$ we have 
 $\lau{SC}{\per{p}}{i}{cL}  =\lau{SC}{\per{p}}{i}{L \times ]0,1[} $ which 
 gives 
 $
 \lau{\IH}{\per{p}}{i}{cL} =
 \lau{\IH}{\per{p}}{i}{L}
 $
 if $i \leq \per{t}(\vartheta)- \per{p}(\vartheta)$.

 For a singular simplex $\sigma \colon \Delta^i \to cL$ we define the {\em 
 cone} $c\sigma \colon \Delta^{i+1} \to cL$ by
 $$
 c\sigma (x_0,\ldots,x_{i+1}) = ( 1-x_{i+1}) \cdot 
\sigma\left(\frac{x_0}{1-x_{i+1}}, \ldots , 
 \frac{x_i}{1-x_{i+1}}\right).
 $$
 Here, we have written $r \cdot [x,s] = [x,rs]$ for a point $[x,s] \in cL$ and a number 
 $r\in [0,1]$.
 In the same way, we define the cone $c\xi$ of a singular chain $\xi$.
It defines the linear operator 
\begin{equation}
    \label{conop}
    c \colon \lau{AC}{\per{p}}{\geq 1+
\per{t}(\vartheta) - \per{p}(\vartheta)}{cL} 
\TO \lau{AC}{\per{p}}{\geq 
2+\per{t}(\vartheta) - \per{p}(\vartheta)}{cL}.
\end{equation}
Let us  prove this property.
Take $ \sigma  \in \lau{AC}{\per{p}}{\geq
1 + \per{t}(\vartheta) - \per{p}(\vartheta)}{cL} $ and prove that
$ c\sigma  \in \lau{AC}{\per{p}}{\geq
2 + \per{t}(\vartheta) - \per{p}(\vartheta)}{cL} $.
 Notice first that, for $x_{i+1}\ne 1$, we have 
 $$
 \underbrace{\left(\frac{x_0}{1-x_{i+1}}, \ldots , 
  \frac{x_i}{1-x_{i+1}}\right)}_{\tau(x_{0},\ldots,x_{i+1}) }
  \in j-\hbox{skeleton of } \Delta^{i} \Longrightarrow
  (x_{0},\ldots,x_{i+1}) 
  \in (j+1)-\hbox{skeleton of } \Delta^{i+1} .
  $$
  So,
  \begin{eqnarray*}
      (c\sigma)^{-1}(\vartheta) & = & \left\{(0,\ldots,0,1)\right\} \cup   
      \left\{(x_0, \ldots, x_{i+1}) \in \Delta^{i+1} \ / \ \tau(x_{0},\ldots,x_{i+1})
      \in \sigma^{-1}(\vartheta)\hbox{ and }x_{i+1}\ne 
  1 \right\}\\[,2cm]
       & \subset & \left\{(0,\ldots,0,1)\right\} \cup   
       (i-2-\per{t}(\vartheta) +\per{p}(\vartheta)+1)-\hbox{skeleton 
       of } \Delta^{i+1}\\[,2cm]
       & \subset &(i+1-2-\per{t}(\vartheta) +\per{p}(\vartheta))-\hbox{skeleton 
       of } \Delta^{i+1},
  \end{eqnarray*}
  since $i\geq 1+\per{t}(\vartheta) - \per{p}(\vartheta)$. For
 each singular stratum $S$ of $L$ we have:
  \begin{eqnarray*}
       (c\sigma)^{-1}(S \times ]0,1[) & = &  
       \left\{(x_0, \ldots, x_{i+1}) \in \Delta^{i+1} \ / \ \tau(x_{0},\ldots,x_{i+1})\in \sigma^{-1}(S \times 
       ]0,1[) \hbox{ and }x_{i+1}\ne 
   1 \right\}\\[,2cm]
	& \subset & (i+1-2-\per{t}(\vartheta) +\per{p}(\vartheta))-\hbox{skeleton 
	of } \Delta^{i+1}.
   \end{eqnarray*}
   We conclude that $ c\sigma  \in \lau{AC}{\per{p}}{\geq
2 + \per{t}(\vartheta) - \per{p}(\vartheta)}{cL} $.
 Notice that any singular simplex verifies the formula
 $
 \partial c \sigma  = 
 c\partial \sigma + (-1)^{i+1}\sigma,
 $
 if $i > 0$ and  $
 \partial c \sigma  = 
\vartheta -\sigma$, if $i =0$.
 
 Consider now a cycle $ \xi \in \lau{SC}{\per{p}}{i}{cL}$ with  $i \geq  
 1+\per{t}(\vartheta)- 
 \per{p}(\vartheta)$ and $i\ne 0$. Since $\xi = (-1)^{i+1}\partial c \xi$ then 
 $c\xi \in \lau{SC}{\per{p}}{i+1}{cL}$  and therefore
 $\lau{\IH}{\per{p}}{i}{cL}=0$. 

 For $i=0\geq  1+ \per{t}(\vartheta)- 
 \per{p}(\vartheta)$ we get 
 that for any point $\sigma$ of $cL - \{ \vartheta \}$ the cone $c\sigma$  is a 
 $\per{p}$-allowable chain with
 $\partial c\sigma = \sigma - \vartheta$. This gives  $\lau{\IH}{\per{p}}{0}{cL}=\R$.
 \qed

 In some cases, the intersection homology can be expressed in terms of the usual 
 homology $\homo{H}{*}{-}$  (see \cite{GM1}).

\bP
\label{part}
Let $X$ be a stratified pseudomanifold. Then 
\begin{itemize}
    \item$\lau{\IH}{\per{p}}{*}{X} = \homo{H}{*}{X-\Sigma_{X}}$ if 
    $\per{p} < \per{0}$, and
    
    \item  $\lau{\IH}{\per{q}}{*}{X} = \homo{H}{*}{X}$ if 
    $\per{q} \geq \per{t}$ and $X$ is normal.
\end{itemize}
\eP
\pro
We prove, by induction on the depth, that the natural inclusions 
$I_{X}\colon \homo{S}{*}{X-\Sigma_{X}} \hookrightarrow \lau{SC}{\per{p}}{*}{X}$ and
$J_{X}\colon \lau{SC}{\per{q}}{*}{X} \hookrightarrow  \homo{S}{*}{X} $ 
are quasi-isomorphisms (i.e. isomorphisms in cohomology). When the 
depth of $X$ is 0 then 
the above inclusions are, in fact, two identities. In the general case, 
we suppose that the result is true for each link $L_{S}$ of $X$ and 
we proceed in two steps.

\medskip

{\em \underline{First Step} The operators $I_{V}$ and $J_{V}$ are 
quasi-isomorphisms when $V$ is an open subset of a chart $(U,\phii)$ of $X$.}

\smallskip

\nt First of all, we identify the open subset $U$ with the product $\R^n \times cL_{S}$ 
through 
$\phii$.
We define a {\em cube}  as a product 
$]a_{1},b_{1}[ \times \cdots \times 
]a_{n},b_{n}[ \subset \R^n$.
The {\em truncated cone} $c_{t}L_{S}$ is  
the quotient $c_{t}L_{S} = L_{S}\times [0,t[ / L_{S}\times \{ 0\}$.
Consider the open covering
$$
\mathcal{V} = \Big\{ C\times c_{t}L_{S} \subset V\ / \ C \hbox{ cube, } 
t \in ]0,1[  \Big\}  \cup \Big\{ C \times L_{s} \times ]a,b[  \subset 
V \ / \ 
C \hbox{ cube, } a,b 
\in ]0,1[  \Big\}
$$
of $V$. Notice that this family is closed for finite intersections. 

We use the 
Bredon's trick relatively to the covering $\mathcal{V}$ and to the 
statement
$$
Q(W) = ``\hbox{ The operators $I_{W}$ and $J_{W}$ are 
quasi-isomorphisms}''
$$
(cf. Lemma \ref{Brt}). Let us verify the properties (BT1), (BT2) and 
(BT3).

\begin{itemize}
    \item[]
    \begin{itemize}
\item[(BT1)] From the product formula and the induction hypothesis, it suffices to prove that the 
operators 
$I_{cL_{S}}$ and $J_{cL_{S}}$ are quasi-isomorphisms. This comes from:
\begin{itemize}
    \item  $\lau{\IH}{\per{p}	}{*}{cL_{S}} \stackrel{\ref{cono1}}{=} 
    \lau{\IH}{\per{p}}{\leq \per{t}(\vartheta) - \per{p}(\vartheta)}{L_{S}}
    \stackrel{\per{p}(\vartheta) < 0}{= \! = \! =} 
    \lau{\IH}{\per{p}}{*}{L_{S}} \stackrel{ind}{=} 
    \homo{H}{*}{L _{S}- \Sigma_{L_{S}}}     \stackrel{prod}{= \! =}  \homo{H}{*}{cL_{S} - \Sigma_{cL_{S}}}$.

    \item  For $\per{q}(\vartheta) = \per{t}(\vartheta) $ we have: 
    $\lau{\IH}{\per{q}}{*}{cL_{S}} 
    \stackrel{\ref{cono1}}{=} 
    \lau{\IH}{\per{q}}{0}{L_{S}}
    \stackrel{ind}{=}
	\homo{H}{0}{L_{S}}
	\stackrel{norm}{=} 
	\R =     \homo{H}{*}{cL_{S} }$.
	
	\item  For $\per{q}(\vartheta) > \per{t}(\vartheta) $ we have: 
	$\lau{\IH}{\per{q}}{*}{cL_{S}} \stackrel{\ref{cono1}}{=}
	    \R 
	    =     \homo{H}{*}{cL_{S} }$.
\end{itemize}
\item[(BT2)] Mayer-Vietoris.
\item[(BT3)] Straightforward.
\end{itemize}
\end{itemize}

\medskip

{\em \underline{Second Step.} The operators $I_{X}$ and $J_{X}$ are 
quasi-isomorphisms.}

\smallskip

\nt Consider the open covering
$
\mathcal{V} = \Big\{ V  \hbox{ open subset of a chart $(U,\phii)$ of $X$}\Big\}
$
of $X$. Notice that this family is closed for finite intersections.
We use the 
Bredon's trick relatively to the covering $\mathcal{V}$ and to the 
statement
$
Q(W) = \hbox{ ``The operators $I_{W}$ and $J_{W}$ are 
quasi-isomorphisms''}
$
(cf. Lemma \ref{Brt}). Let us verify the properties (BT1), (BT2) and 
(BT3).

\begin{itemize}
    \item[]
    \begin{itemize}
\item[(BT1)] First Step.
\item[(BT2)] Mayer-Vietoris.
\item[(BT3)] Straightforward. \qed
\end{itemize}
\end{itemize}

\prgg {\bf Remark}. Notice that we can replace the normality of $X$ 
by the connectedness of the links $\{ L_{S }\ / \ \per{q}(S) = \per{t}(S)\}$.

\medskip

The following result will be needed in the last section.

\bC
\label{norm}
Let $X$ be a connected normal stratified pseudomanifold. Then,  for any
perversity $\per{p}$,
 we have $\lau{\IH}{\per{p}}{0}{X} = \R.$
\eC
\pro
We prove this result by induction on the depth. When $\depth X =  0$
then $\lau{\IH}{0}{\per{p}}{X}\stackrel{\Sigma_{X = 
\emptyset}}{= \! = \! =} \homo{H}{0}{X} = \R$.
Consider now the general case.
Notice that any point $\sigma \in X - \Sigma_{X}$  is a 
$\per{p}$-intersection cycle. So $\lau{\IH}{\per{p}}{0}{X}\ne 0$. 
We prove that $[\sigma_{1}] = [\sigma_{2}]$ in $\lau{\IH}{\per{p}}{0}{X}$ 
for two $\per{p}$-allowable points. This is the case when 
$\sigma_{1}, \sigma_{2} \in X - 
\Sigma_{X}$, since we know from \cite{P} that $X-\Sigma_{X}$ is
connected. Consider now a $\per{p}$-intersection 
cycle $\sigma_{1} \in \Sigma_{X}$ and $\phii \colon U  \to \R^n\times cL_S$ 
a chart of $X$ containing $\sigma_{1}$.
Since
$$
\lau{\IH}{\per{p}}{0}{U} = \lau{\IH}{\per{p}}{0}{\R^n\times cL_S}
\stackrel{prod}{= \! =} \lau{\IH}{\per{p}}{0}{cL_S} 
\stackrel{\ref{cono1}, \ ind}{= \! = \! = \! = \! =} \R,
$$
we have $[\sigma_1] = [\sigma_2]$ in $\lau{\IH}{\per{p}}{0}{X}$ for 
some $\per{p}$-intersection point $\sigma_{2} \in X - \Sigma_{X}$. 
Then $\lau{\IH}{\per{p}}{0}{X} = \R$.
\qed
\prgg {\bf Relative case}. The conical formula given by 
Proposition \ref{cono1} for the intersection homology differs from that of Proposition 
\ref{cono5} for the intersection cohomology: we do not have 
$\lau{\IH}{*}{\per{p}}{cL} =  \lau{\IH}{\per{t}-\per{p}}{*}{cL} $ when the 
perversity $\per{p}$ is not positive. It is natural to think that the 
closed saturated subset
$$
\ib{X}{\per{p}} = \bigcup_{ S \preceq S_{1} \atop
\per{p} (S_{1})< 0}S 
=
\bigcup_{ \per{p} (S)< 0}\overline{S}  \stackrel{loc \ finit}{= \! = \! =\! = \! = } 
\overline{\bigcup_{ \per{p} 
(S)< 0}S}
$$
plays a key rôle in the de Rham Theorem. This is indeed the case.

The subset $\ib{X}{\per{p}}$ is a stratified pseudomanifold where the maximal strata 
may have different dimensions. For any perversity $\per{q}$ (on $X$) 
we have the notion of $\per{q}$-allowable chain as in 2.2. We denote by 
$\lau{AC}{\per{q}}{*}{\ib{X}{\per{p}}}$ the complex of these $\per{q}$-allowable 
chains. Equivalently, 
$$
\lau{AC}{\per{q}}{*}{\ib{X}{\per{p}}} =
\homo{S}{*}{\ib{X}{\per{p}}} \cap 
\lau{AC}{\per{q}}{*}{X}.
$$

In order to recover the de Rham Theorem we introduce 
the following notion of  relative intersection homology.
We denote by
$
\lau{SC}{\bar q}{*}{X,\ib{X}{\per{p}}}
$
the differential complex 
$$
\fracc{\left(\lau{AC}{\per{q}}{*}{X}+\lau{AC}{\per{q}+ 
\per{1}}{*}{\ib{X}{\per{p}} }\right)
\cap \partial^{-1}\left(\lau{AC}{\per{q}}{*-1}{X} 
+\lau{AC}{\per{q}+ \per{1}}{*-1}{\ib{X}{\per{p}} }\right)}
{\lau{AC}{\per{q}+ \per{1}}{*}{\ib{X}{\per{p}} } \cap \partial^{-1}
\left(\lau{AC}{\per{q}+ \per{1}}{*-1}{\ib{X}{\per{p}}}\right)}$$
and by
$\lau{\IH}{\per{q}}{*}{X,\ib{X}{\per{p}} }$ its cohomology. Of course, 
we have
$\lau{\IH}{\per{q}}{*}{X,\ib{X}{\per{p}} } = \lau{\IH}{\per{q}}{*}{X}$ 
when $\ib{X}{\per{p}} = 
\emptyset$ and 
$\lau{\IH}{\per{q}}{*}{X,\ib{X}{\per{p}} } = \homo{H}{*}{X,\ib{X}{\per{p}} }$ when 
$\per{q} \geq \per{t}+\per{2}$ (see also \refp{ZZZ}).

Since the complexes defining the relative complex $\lau{SC}{\per{q}}{*}{X,\ib{X}{\per{p}}}$ 
verify the Mayer-Vietoris formula then the relative cohomology also 
verifies this property.
For the same reason we have the product formula 
$\lau{\IH}{\per{p}}{*}{\R \times X, \R \times \ib{X}{\per{p}}}  = 
\lau{\IH}{\per{p}}{*}{X,\ib{X}{\per{p}}} $.
For the typical local calculation we have the following result. 

\bC
\label{cono2}
Let $L$ be a compact stratified pseudomanifold. Then, 
for any perversity $\per{p}$, we have:
\begin{equation}
    \label{equa}
    \lau{\IH}{\per{t} - \per{p}}{i}{cL,\ib{(cL)}{\per{p}} } =
\left\{\begin{array}{cl}
\lau{\IH}{\per{t} -\per{p}}{i}{L,\ib{L}{\per{p}} } & \hbox{if $i \leq \per{p}(\vartheta)$} 
\\[,2cm]
0& \hbox{if $i \geq 1 +  \per{p}(\vartheta)$}.
\end{array}
\right.
\end{equation}
\eC
\pro When $\per{p} \geq \per{0}$ then $\ib{(cL)}{\per{p}}  = 
   \ib{L}{\per{p}}  = \emptyset$  and \refp{equa} comes directly from 
   Lemma \ref{cono1}. Let us suppose $\per{p} \not\geq \per{0}$, 
   which gives $\ib{(cL)}{\per{p}}  = 
   c\left( \ib{L}{\per{p}} \right) \ne\emptyset$, with $c\emptyset = \{ \vartheta 
   \}$. We also use the following equalities:
  \begin{equation}
      \label{cLL}
  \lau{AC}{\per{t} -\per{p}}{j}{cL} = \lau{AC}{\per{t} -\per{p}}{j}{L\times ]0,1[}
  \ \ \hbox{ for }     j \leq \per{p}(\vartheta) + 1,
      \end{equation}
      and
      \begin{equation}
	  \label{cLLL}
	  \lau{AC}{\per{t} -\per{p} + \per{1}}{j}{\ib{(cL)}{\per{p}}} = 
	  \lau{AC}{\per{t} -\per{p}+ \per{1}}{j}{\ib{L}{\per{p}}\times 
	  ]0,1[}
      \ \ \hbox{ for }     j \leq \per{p}(\vartheta),
	  \end{equation}    
   We proceed in four steps following the value of $ i \in 
\N$.

\medskip

{\em \underline{First Step}: 
$i \leq \per{p}(\vartheta)  -1$.}
We have
 $\lau{SC}{\per{t} 
-\per{p}}{j}{cL,\ib{(cL)}{\per{p}} }=\lau{SC}{\per{t} 
-\per{p}}{j}{L \times ]0,1[,\ib{L}{\per{p}} \times ]0,1[ }$ for each  $j \leq \per{p}(\vartheta)  $ (cf. 
\refp{cLL} and \refp{cLLL}) and 
therefore
$\lau{\IH}{\per{t} - \per{p}}{i}{cL,\ib{(cL)}{\per{p}} } =
\lau{\IH}{\per{t} -\per{p}}{i}{L,\ib{L}{\per{p}} }.$

\smallskip

{\em \underline{Second Step}: 
$i = \per{p}(\vartheta)  $.}
The inclusion 
$
\lau{SC}{\per{t} -\per{p}}{*}{L \times ]0,1[,\ib{L}{\per{p}} \times ]0,1[ }
\hookrightarrow \lau{SC}{\per{t} -\per{p}}{*}{cL,\ib{(cL)}{\per{p}}}
$
induces the epimorphism
$
I \colon 
\lau{\IH}{\per{t} -\per{p}}{\per{p}(\vartheta)  
}{L \times ]0,1[,\ib{L}{\per{p}} \times ]0,1[ }
\TO \lau{\IH}{\per{t} -\per{p}}{\per{p}(\vartheta)  }{cL,\ib{(cL)}{\per{p}}}
$
(cf. \refp{cLL} and \refp{cLLL}). It remains to prove that $I$ is a 
monomorphism. Consider $\Big[\overline{\alpha + \beta}\Big] \in \lau{\IH}{\per{t} -\per{p}}{\per{p}(\vartheta)  
}{L \times ]0,1[,\ib{L}{\per{p}} \times ]0,1[ }$ with $I\left( 
\Big[\overline{\alpha + \beta}\Big]\right) =0$. So, 

\Zati $\alpha \in \lau{AC}{\per{t} -\per{p}}{\per{p}(\vartheta)  }{L \times 
    ]0,1[ } \subset \lau{AC}{\per{t} -\per{p} + \per{1}}{\per{p}(\vartheta)  }{L \times 
    ]0,1[ }$,
    
\zati  $\beta \in \lau{AC}{\per{t} 
-\per{p} + \per{1}}{\per{p}(\vartheta)  }{\ib{L}{\per{p}} \times ]0,1[ }
    $,

\medskip 

\nt and there exist

\zati $A\in \lau{AC}{\per{t} -\per{p}}{\per{p}(\vartheta) +1}{cL}
 \stackrel{\refp{cLL} }{=}
 \lau{AC}{\per{t} -\per{p}}{\per{p}(\vartheta) +1}{L 
 \times ]0,1[}$,
    
  \zati $B\in \lau{AC}{\per{t} 
-\per{p} + \per{1}}{\per{p}(\vartheta)   +1}{\ib{\left(cL\right)}{\per{p}} 
}$,

\zati $C\in 
\lau{AC}{\per{t} 
-\per{p}+ \per{1}}{\per{p}(\vartheta)  }{\ib{\left(cL\right)}{\per{p}} }   \cap   \partial^{-1}
  \left(\lau{AC}{\per{t} 
-\per{p}+ \per{1}}{\per{p}(\vartheta)-1}{\ib{\left(cL\right)}{\per{p}} }\right)
  \stackrel{\refp{cLLL} }{=} $ 
\begin{flushright}$\lau{AC}{\per{t} 
-\per{p}+ \per{1}}{\per{p}(\vartheta)  }{\ib{L}{\per{p}} \times ]0,1[ }   \cap   \partial^{-1}
  \left(\lau{AC}{\per{t} 
-\per{p}+ \per{1}}{\per{p}(\vartheta)-1}{\ib{L}{\per{p}} 
\times ]0,1[ }\right)$
\end{flushright}

\medskip 

\nt with 

%
  \zati 
$\alpha + \beta = \partial A + \partial B +C$.

\medskip

\nt Since $\partial A \in \lau{AC}{\per{t} -\per{p} + \per{1}}{\per{p}(\vartheta)  }{L \times 
    ]0,1[ }$ (cf. (c)) then the conditions (a), (b), (d), (e) and (f) give 
    
    \zati
    $\partial B \in \lau{AC}{\per{t} -\per{p} + \per{1}}{\per{p}(\vartheta)  
    }{\ib{L}{\per{p}}  \times 
	]0,1[ }$.
	
	\medskip
	
    \nt We conclude that
    $$
    \partial A = \alpha + (\beta  - \partial B - C )\in 
    \lau{AC}{\per{t} -\per{p} }{\per{p}(\vartheta)  }{L \times 
	]0,1[ }
	+
	\lau{AC}{\per{t} -\per{p} + \per{1}}{\per{p}(\vartheta)  
	}{\ib{L}{\per{p}} \times ]0,1[ },
	$$
	which defines the element 
	$\overline{A} \in 
	\lau{SC}{\per{t} -\per{p} }{\per{p}(\vartheta)  +1 }{L \times 
	]0,1[ ,\ib{L}{\per{p}}  \times 
	]0,1[ }$. If we write $\overline{\partial}$ the derivative of
	$\lau{SC}{\per{t} -\per{p} }{*}{L \times 
	]0,1[ ,\ib{L}{\per{p}}  \times 
	]0,1[ }$ we can write:
	$$
	\Big[\overline{\alpha + \beta}\Big] = 
	\Big[\overline{\partial A + \partial B +C}\Big] 
	=
	\Big[\overline{\partial} \ \overline{A} + \overline{\partial B + C}\Big] 
	\stackrel{(e),(g)}{= \! = \! =}
	\Big[\overline{\partial} \ \overline{A}\Big] 
	= 0.
	$$
Then the operator $I$ is a monomorphism.

   \smallskip 
   
{\em \underline{Third Step}:  $i \geq  1+ \per{p}(\vartheta)$ and $i \ne 0$.}
Consider  a cycle $ \overline{\xi} \in \lau{SC}{\per{t} -\per{p}}{i}{cL,\ib{(cL)}{\per{p}} }$.
Let $\xi = \xi_1 + \xi_2 \in 
\lau{AC}{\per{t} -\per{p}}{i}{cL} + 
\lau{AC}{\per{t} -\per{p}+ \per{1}}{i}{   \ib{(cL)}{\per{p}}  }$ with $\partial \xi \in 
\lau{AC}{\per{t} -\per{p}+ \per{1}}{i-1}{   \ib{(cL)}{\per{p}} }$. Then 
we have 
$c\xi = c\xi_1 + c\xi_2 \in 
\lau{AC}{\per{t} -\per{p}}{i+1}{cL} + 
\lau{AC}{\per{t} -\per{p}+ \per{1}}{i+1}{   \ib{(cL)}{\per{p}}  }$ 
and
$c\partial \xi \in \lau{AC}{\per{t} -\per{p}+ 
\per{1}}{i}{    \ib{(cL)}{\per{p}}}
$
(cf. \refp{conop}). 
Since 
\begin{equation}
    \label{magia}
    \partial \, c \xi = (-1)^{i+1}\xi  + c\partial \xi 
    \end{equation}
    is an element of $ \lau{AC}{\per{t} -\per{p}}{i}{cL}
+\lau{AC}{\per{t} -\per{p}+ 
\per{1}}{i}{    \ib{(cL)}{\per{p}}  }$  then
$\overline{c \xi} \in \lau{SC}{\per{t} -\per{p}}{i+1}{cL,\ib{(cL)}{\per{p}} }$.
This formula gives $\partial c\partial \xi  = (-1)^i \partial \xi$ 
and therefore 
$c\partial \xi \in \lau{AC}{\per{t} -\per{p}+ 
\per{1}}{i}{    \ib{(cL)}{\per{p}}  } \cap \partial^{-1}\left( 
\lau{AC}{\per{t} -\per{p}+ 
\per{1}}{i-1}{    \ib{(cL)}{\per{p}} } 
\right).
$
Applying \refp{magia}  we obtain 
$\left[\ \overline{\xi}\ \right] = 0 $ on 
$\lau{\IH}{\per{t} -\per{p}}{i}{cL,  \ib{(cL)}{\per{p}}}=0$.  

\smallskip 

{\em \underline{Fourth Step}:  $i=0\geq  1+
\per{p}(\vartheta)$}.  
For any point $\sigma  \in \lau{AC}{\per{t} -\per{p}}{0}{ cL }$ the cone 
$c\sigma$  is a 
$(\per{t} -\per{p})$-allowable chain with
$\partial c\sigma = \sigma - \vartheta$. Since the point $\vartheta $
belongs to the complex $\lau{AC}{\per{t} -\per{p}+ 
\per{1}}{0}{ \ib{(cL)}{\per{p}}  }
$
then 
$\lau{\IH}{\per{t} -\per{p}}{0}{cL,\ib{(cL)}{\per{p}} }=0$.
\qed

\prg{\bf Intersection homology. Second approach} (see \cite{S}).  In order to integrate 
differential forms on allowable simplices, we need to introduce some 
amount of
smoothness on these simplices. Since $X$ is not a manifold, we work 
in the manifold $\wt{X}$. In fact we consider those allowable simplices which 
are liftable to  smooth simplices in
$\wt{X}$.

\prgg {\bf Linear unfolding.} The {\it unfolding} of the standard
simplex $\Delta$, relative to the decomposition
$\Delta = \Delta_0 * \cdots * \Delta_j$, is the map  
$
\mu_\Delta\colon \wt{\Delta} = \bar{c}\Delta_0
\times
\cdots
\times \bar{c}\Delta_{j-1} \times  \Delta_j
\longrightarrow
\Delta $  defined by  
\begin{eqnarray*}
 \mu_\Delta ([x_0,t_0], \ldots, [x_{j-1},t_{j-1}],x_j) &= &
t_0x_0+(1-t_0)t_1x_1+\cdots+(1-t_0)\cdots(1-t_{j-2})t_{j-1}x_{j-1}\\[,2cm]
&&+(1-t_0)\cdots(1-t_{j-1})x_j,
\end{eqnarray*}
where
$\bar{c}\Delta_i$ denotes the closed cone  
$\Delta_i \times [0,1] \big/
\Delta_i \times \{ 0\}$ and $[x_i,t_i]$ a point of it.  This map is smooth and its
restriction $\mu_\Delta\colon
\inte(\wt{\Delta})  \longrightarrow \inte(\Delta) $ is a diffeomorphism ($\inte(P) = P-\partial P$ is the {\em interior} of the polyhedron
$P$). It sends a face $U$  of
$\wt{\Delta}$ to a face $V$ of $\Delta$ and the restriction
$\mu_\Delta \colon \inte(U) \to \inte(V)$ is a submersion.

On the boundary  $\partial\wt{\Delta}$ we find not only the blow-up
$\wt{\partial\Delta}$ of the boundary $
\partial\Delta$ of $\Delta$ but also  the faces  
$$F =
\bar{c}\Delta_0 \times  \cdots \times
\bar{c}\Delta_{i-1} \times ( \Delta_i
\times \{ 1\}) \times \bar{c}\Delta_{i+1} \times 
\cdots \times
\bar{c}\Delta_{j-1} \times \Delta_j$$ with  $i\in \{ 0,\ldots,j-2\}$ or $i=j-1$ and
$\dim \Delta_j >0$, which we  call {\em bad faces}.  This gives 
the decomposition
\begin{equation}
    \partial \wt\Delta = \wt{\partial \Delta} + \delta \wt\Delta
    \label{descomp}
\end{equation}
Notice that

\begin{equation}
    \label{bad}\dim \mu_\Delta ( F  )= \dim( \Delta_0 * \cdots *
\Delta_i)< \dim \Delta -1 =\dim F.
\end{equation}

\prgg{\bf Liftable simplices.}
A {\it liftable simplex} is a singular
simplex $\sigma \colon
\Delta \to X$ verifying the  following two conditions. 
\begin{itemize}
\item[]
\begin{itemize}
\item[(Lif1)] {\em Each pull back $\sigma^{-1}(X_i)$ is a face of $\Delta$}.  

\item[(Lif2)]  {\em There exists a decomposition 
$\Delta = \Delta_0 * \cdots * \Delta_j$ and a smooth map
(called {\em lifting})
$\wt{\sigma}
\colon \wt{\Delta} \to \wt{X}$  with ${\cal L}_X \rond
\wt{\sigma }=
\sigma \rond\mu_\Delta$}.
\end{itemize}
\end{itemize}
A singular chain $\xi = \sumas{j=1}{m} r_{j} 
\sigma_{j}$ is {\em liftable} if each singular simplex $\sigma_j$ is 
liftable. Since a face of  a liftable simplex is again a liftable 
simplex then the family
$\homo{L}{*}{X}$ of liftable chains is a differential complex. We 
denote by
$
\lau{LC}{\per{p}}{*}{X} =
\lau{AC}{\per{p}}{*}{X} \cap \homo{L}{*}{X}
$
the graded vector space of the $\per{p}$-allowable liftable chains and 
by
$
\lau{RC}{\per{p}}{*}{X} = \lau{LC}{\per{p}}{*}{X} \cap \partial^{-1} 
\lau{LC}{\per{p}}{*-1}{X}
$
the associated 
 differential complex.
The 
cohomology of this complex verifies  the Mayer-Vietoris formula and the product formula 
$\homo{H}{*}{\lau{RC}{\per{p}}{\cdot}{\R \times X} } = 
\homo{H}{*}{\lau{RC}{\per{p}}{\cdot}{X} }$.
For the typical local calculation we have the following result. It corrects  Proposition 2.2.5 of \cite{S}.
\bP
\label{cono3}
Let $L$ be a compact unfoldable pseudomanifold.
Consider on $cL$ the canonical induced unfolding.  Then
$$\homo{H}{i}{\lau{RC}{*}{\bar p}{cL}}=
\left\{\begin{array}{cl}
\homo{H}{i}{\lau{RC}{*}{\bar p}{L}}& \hbox{if $i \leq \per{t}(\vartheta)- \per{p}(\vartheta)$} 
\\[,2cm]
0& \hbox{if $0 \ne i \geq 1+\per{t}(\vartheta)- \per{p}(\vartheta)$}
\\[,2cm]
\R & \hbox{if $0 = i \geq 1+\per{t}(\vartheta)- \per{p}(\vartheta)$}.
\end{array}
\right.
$$
\eP
\pro
We proceed as in Proposition \ref{cono1}. In fact, it suffices to prove that the 
cone $c\sigma \colon \overline{c}\Delta \to cL$ of 
a $\per{p}$-allowable liftable simplex $\sigma \colon \Delta \to cL$, with 
$\dim \Delta \geq 1+\per{t}(\vartheta) -\per{p}(\vartheta)$, 
is a $\per{p}$-allowable liftable 
simplex. Let us verify properties (All)$_{c\sigma}$, (Lif1)$_{c\sigma}$ and (Lif2)$_{c\sigma}$.

Put $\overline{c}\Delta = \Delta * \{ Q \}$ with $c\sigma (tP + (1-t)Q)=
t  \cdot \sigma(P) $. We have:
$$
(c\sigma)^{-1}\left(cL\right)_i =
\left\{
\begin{array}{ll}
\{ Q\} & \hbox{if } \sigma^{-1}(cL)_{i} = \emptyset\\[,2cm]
\overline{c}(\sigma^{-1}(cL)_{i}) & \hbox{if } \sigma^{-1}(cL)_{i} \not= \emptyset,
\end{array}
\right.
$$
for $i \geq 0.$

We obtain (Lif1)$_{c\sigma}$ from (Lif1)$_{\sigma}$. 
To prove the property (All)$_{c\sigma}$ we  consider a 
stratum $S \in \mathcal{S}_{cL}$. We have
\begin{eqnarray*}
(c\sigma)^{-1}(S) &=&
\left\{
\begin{array}{ll}
\{ Q\} & \hbox{if } S = \{ \vartheta\} ;\sigma^{-1}(\vartheta ) = \emptyset \\[,2cm]
\overline{c}\left(\sigma^{-1}( \vartheta) \right) &\hbox{if } S = \{ \vartheta\}; 
\sigma^{-1}(\vartheta ) \ne \emptyset \\[,2cm]
\emptyset
&\hbox{if } S \ne \{ \vartheta\}; \sigma^{-1}(S) =\emptyset  \\[,2cm]
\overline{c}\left(\sigma^{-1}( \vartheta) \right) - \{ Q \} 
&\hbox{if } S \ne \{ \vartheta\} ; \sigma^{-1}(S) \ne \emptyset \\[,2cm]

\end{array}
\right.
\\
&\stackrel{\hbox{(All)}_{\sigma}}{\subset} &
\left\{
\begin{array}{l}
0-\hbox{skeleton of } \overline{c}\Delta 
\\[,2cm]
(1+(\dim \Delta - 2 - \per{t}(\vartheta)+ \per{p}(\vartheta) ))-\hbox{skeleton of } \overline{c}
\Delta
\\[,2cm]
\emptyset \\[,2cm]
(1+(\dim \Delta - 2 - \per{t}(S)+ \per{p}(S) ))-\hbox{skeleton of } \overline{c}\Delta
\end{array}
\right.
\\[,2cm]
&\subset &
(\dim \overline{c}\Delta - 2 - \per{t}(\vartheta)+ 
\per{p}(\vartheta) ))-\hbox{skeleton of } \overline{c}\Delta,
\end{eqnarray*}
since $\dim \Delta \geq 1+\per{t}(\vartheta) - \per{p}(\vartheta)$.
Now we prove (Lif2)$_{c\sigma}$. 
Consider the decomposition  $\Delta = \Delta_0 * \cdots *
\Delta_j $ given by $\sigma$, and the 
smooth map $\wt{\sigma} = (\wt{\sigma} _1,\wt{\sigma} _2) \colon 
\wt{\Delta} \to \wt{L} \times ]-1,1[$ given by (Lif2)$_\sigma$. We 
have the decomposition  $\overline{c}\Delta =  \{Q\} * \Delta_0 * \cdots *
\Delta_j $ 
whose lifting 
$
\mu_{\overline{c}\Delta} \colon \wt{\overline{c}\Delta} = \overline{c}\{Q\} \times \wt{\Delta} \TO 
\overline{c}\Delta
$
is defined by
$$
\mu_{\overline{c}\Delta}  ([Q,t],x) = tQ +(1-t) \mu_\Delta (x).
$$
Let $\wt{c\sigma} \colon  \overline{c}\{Q\}\times \wt{\Delta}  \TO \wt{L} 
\times ]-1,1[$ be  the smooth map defined by 
$$
\wt{c\sigma}( [Q,t],x)  = (\wt{\sigma} _1(x), (1-t)  \cdot \wt{\sigma} _2(x)).
$$
Finally, for each $([Q,t],x) \in \overline{c}\{Q\}\times \wt{\Delta}$ we have:
\begin{eqnarray*}
c\sigma \mu_{\overline{c}\Delta}([Q,t],x)  &= &c\sigma(tQ +(1-t) \mu_\Delta (x))
= (1-t) \cdot \sigma \mu_\Delta (x) = (1-t) \cdot \mathcal{L}_{cL} \wt{\sigma}(x)
\\[,2cm]
&=&(1-t)\left[ \mathcal{L}_{L} \wt{\sigma}_1(x) ,|\wt{\sigma}_2(x)|
\right]=
\mathcal{L}_{cL} \wt{c\sigma}([Q,t],x).
\end{eqnarray*}
This gives  (Lif2)$_{c\sigma}$.  \qed

\prgg{\bf Relative case.} Following 2.2.5 we consider
$$
\lau{LC}{\per{q}}{*}{\ib{X}{\per{p}}} =
\homo{S}{*}{\ib{X}{\per{p}}} \cap 
\lau{LC}{\per{q}}{*}{X}.
$$
and we define the relative complex
$
\lau{RC}{\bar q}{*}{X,\ib{X}{\per{p}}}
$
by
$$
\fracc{\left(\lau{LC}{\per{q}}{*}{X}+\lau{LC}{\per{q}+ 
\per{1}}{*}{\ib{X}{\per{p}} }\right)
\cap \partial^{-1}\left(\lau{LC}{\per{q}}{*-1}{X} 
+\lau{LC}{\per{q}+ \per{1}}{*-1}{\ib{X}{\per{p}} }\right)}{\lau{LC}{\per{q}+ 
\per{1}}{*}{\ib{X}{\per{p}} } \cap \partial^{-1} \left( 
\lau{LC}{\per{q}+ 
\per{1}}{*-1}{ \ib{X}{\per{p}} } \right)}.
$$
We have
$\lau{LC}{\per{q}}{*}{X,\ib{X}{\per{p}} } = \lau{LC}{\per{q}}{*}{X}$ 
when $\ib{X}{\per{p}} = 
\emptyset$.
Since the complexes defining the relative complex $\lau{RC}{\bar p}{*}{X,Z}$ 
verify the Mayer-Vietoris formula then the relative cohomology also 
verifies this property.
We have, for the same reason, the product formula 
$
\homo{H}{*}{\lau{RC}{\per{p}}{\cdot}{\R \times X, \R \times Z} } = 
\homo{H}{*}{\lau{RC}{\per{p}}{\cdot}{X,  Z} } $.
For the typical local calculation we have (see \cite{S}): 
\bC
\label{cono4}
Let $L$ be a compact unfoldable pseudomanifold. 
Consider on $cL$ the canonical induced unfolding. Then
$$
\homo{H}{i}{\lau{RC}{\per{t} - \per{p}}{*}{cL,\ib{(cL)}{\per{p}} }}=
\left\{\begin{array}{cl}
\homo{H}{i}{\lau{RC}{\per{t} - \per{p}}{*}{L,\ib{L}{\per{p}}}}& \hbox{if $i  \leq 
\per{p}( \vartheta)$} 
\\[,2cm]
0& \hbox{if $i \geq 1 + \per{p}(\vartheta)$}.
\end{array}
\right.
$$
\eC
\pro The same proof as that of Corollary \ref{cono2}.
\qed

\prg {\bf Comparing the two approaches}.
When the perversity $\per{p}$ lies between $\per{0} $ and $\per{t}$ then we 
have the isomorphism $\lau{\IH}{\per{p}}{*}{X} 
= \homo{H}{*}{\lau{RC}{\per{p}}{\cdot}{X}}$ (cf. \cite{S}).
We are going to check that this property extends to any perversity for the 
absolute and the relative case.
\bP
\label{iso1}
For any perversity $\per{p}$, the inclusion $\lau{RC}{\per{p}}{*}{X} 
\hookrightarrow 
\lau{SC}{\per{p}}{*}{X}$ induces the isomorphism 
$
\homo{H}{*}{\lau{RC}{\per{p}}{\cdot}{X}}=\lau{\IH}{\per{p}}{*}{X} $.
\eP
\pro
We proceed by induction on the depth. When $\depth X = 0$ then 
$\lau{SC}{\per{p}}{*}{X} = 
\lau{RC}{\per{p}}{*}{X} = \homo{S}{*}{X} $. In the general case, 
we use the Bredon's trick (see the proof of Proposition \ref{part}) and we reduce the problem to a chart 
$X = \R^n \times cL_{S}$. Here, we apply the  product formula and we 
reduce the problem
to $X = cL_{S}$. We end the proof by applying Propositions  \ref{cono1},
 \ref{cono3}  and the induction hypothesis.
\qed

\bP
\label{iso2}
 For any perversity $\per{p}$, 
the inclusion $\lau{RC}{\per{t} -\per{p}}{*}{X,\ib{X}{\per{p}}} \hookrightarrow 
\lau{SC}{\per{t} -\per{p}}{*}{X,\ib{X}{\per{p}}}$ induces the isomorphism 
$\homo{H}{*}{\lau{RC}{\per{t} -\per{p}}{\cdot}{X,\ib{X}{\per{p}}}}=
\lau{\IH}{\per{t} -\per{p}}{*}{X,\ib{X}{\per{p}}} $.
\eP
\pro
We proceed by induction on the depth. If $\depth X = 0$ then 
$\lau{SC}{\per{t} -\per{p}}{*}{X,\ib{X}{\per{p}}} \! = \! 
\lau{RC}{\per{t} -\per{p}}{*}{X,\ib{X}{\per{p}}}$ $=$ 
$\homo{S}{*}{X} $. In the general case, 
we use the Bredon's trick (see the proof of Proposition \ref{part}) and we reduce the problem to a chart 
$X = \R^n \times cL_{S}$ with $\ib{X}{\per{p}}=  \R^n \times 
\ib{(cL_{S})}{\per{p}}$.
Here, we apply the  product formula and we reduce the problem
to 
$(X ,\ib{X}{\per{p}}) = \left( cL_{S},\ib{(cL_{S})}{\per{p}}\right)$.
We end the proof by applying 
 Corollaries  \ref{cono2}, 
\ref{cono4}  and the induction hypothesis.
\qed

\section{Intersection cohomology}
The de Rham intersection cohomology was introduced by Brylinski in 
\cite{Bry}. In our paper we use the presentation of \cite{S}. 

\prg {\bf Perverse forms.} A {\em liftable form} is a differential 
form $\om \in \coho{\Om}{*}{X-\Sigma_{X}}$ possessing a 
{\em lifting}, that is, a 
differential form $\wt\om \in \coho{\Om}{*}{\wt{X}}$
verifying $\wt\om = 
\mathcal{L}_{X}^{*}\om$ on $\mathcal{L}^{-1}_{X}(X - \Sigma_{X})$.

Given two liftable differential forms $\om$, $\eta$ we have the equalities:
\begin{equation}
    \label{elem}
    \wt{\om + \eta} = \wt\om + \wt\eta \qquad ; \qquad
    \wt{\om \wedge \eta} = \wt\om \wedge \wt\eta 
    \qquad ; \qquad \wt{d\om} = d \wt\om.
\end{equation}
We denote by $\coho{\Pi}{*}{X}$ the differential complex of liftable forms.

Recall that,  for each singular 
stratum $S$, the restriction $\mathcal{L}_{S} \colon 
\mathcal{L}_{S}^{-1}(S) \TO S$ is a fiber bundle. 
For a differential 
form $\eta \in \coho{\Om}{*}{\mathcal{L}_{S}^{-1}(S) }$ we define its {\em vertical 
degree}
as
$$
v_{S}( \eta ) = \min \left\{ j \in \N \ \Big/ \  \begin{array}{l}
i_{\xi_{0} }\cdots i_{\xi_{j} }\eta =  0 \hbox{ for 
each family of vector fields}  \\ \xi_{0}, \ldots \xi_{j}
\hbox{ tangent to 
fibers of }\mathcal{L}_{S} \colon 
\mathcal{L}_{S}^{-1}(S) \TO S \end{array} \right\}
$$
(cf. \cite{Bry},\cite{S}).  
The {\em perverse degree} 
$\| \om \|_{S}$ of $\om$ relative to $S$ is 
the vertical degree of the restriction $\wt\om$ relatively to $\mathcal{L}_{S} \colon 
\mathcal{L}_{S}^{-1}(S) \TO S$, that is,  $$\| \om \|_{S} = 
v_{S}\left( \wt\om|_{\mathcal{L}_{S}^{-1}(S)}\right).$$
The differential complex of 
$\per{p}$-intersection differential forms is 
$$\lau{\Om}{*}{\per{p}}{X} = 
\{
\om \in \coho{\Pi}{*}{X} \ / \max 
\left( ||\om||_{S}, ||d\om||_{S}\right) \leq  \per{p}(S) \ \forall 
\hbox{ singular stratum } S
\}.
$$
The cohomology $\lau{\IH}{*}{\per{p}}{X}$ of this complex is the 
{\em $\per{p}$-intersection 
cohomology } of $X$.
The intersection cohomology verifies two important 
computational properties: the Mayer-Vietoris property and the product formula 
$\lau{\IH}{*}{\per{p}}{\R \times X}  = \lau{\IH}{*}{\per{p}}{X}$. 
The usual local calculations   (see 
\cite{Bry},
\cite{S}) give
\bP
\label{cono5}
Let $L$ be a compact stratified pseudomanifold. Then
$$\lau{\IH}{i}{\per{p}}{cL} =
\left\{\begin{array}{ll}
\lau{\IH}{i}{\per{p}}{L} & \hbox{if $i \leq  \per{p}(\vartheta)$} 
\\[,2cm]
0& \hbox{if $i > \per{p}(\vartheta)$}.
\end{array}
\right.
$$
\eP
\prg {\bf Integration}. The relationship between the intersection 
homology and cohomology is established by using the integration of 
differential forms on simplices. Since $X$ is not a manifold, we work 
on the blow up $\wt{X}$.

Consider a $\phii \colon \Delta \to X$ 
a liftable simplex. We know that there exists a stratum $S$ 
containing $\sigma(\inte (\Delta))$. Since $\mu_{\Delta }\colon 
\inte(\wt{\Delta}) \to 
\inte (\Delta)$ is a diffeomorphism then $\sigma = \mathcal{L}_{X} 
\rondp \wt\sigma \rondp \mu_{\Delta}^{-1} \colon \inte ( \Delta) 
\TO S$ is a smooth map.

Consider now a liftable differential form $\om \in 
\coho{\Pi}{*}{X}$ and define the {\em integration} as
\begin{equation}
    \label{def}
    \int_{\sigma}\om = \left\{
    \begin{array}{ll}
        {\displaystyle \int_{\inte(\Delta)}\sigma^* \om}& 
        \hbox{if $S$ a regular stratum (i.e. $\sigma(\Delta) 
        \not\subset \Sigma_{X}$)} \\[,6cm]
        0 & \hbox{if $S$ a singular stratum (i.e. $\sigma(\Delta) 
       \subset \Sigma_{X}$)}
    \end{array}
    \right.
\end{equation}
This definition makes sense since
\begin{equation}
    \label{Jos}
\int_{\inte(\Delta)}\sigma^* \om
=
\int_{\inte(\wt{\Delta})} {\wt\sigma}^{*}\wt\om
=
\int_{\wt{\Delta}} {\wt\sigma}^{*}\wt\om.
\end{equation}
By linearity, we have the linear pairing
$
    {\displaystyle \int \colon \coho{\Pi}{*}{X} \TO \Hom 
    \left( 
    \homo{L}{*}{X},\R \right)}.
$
This operator commutes with the differential $d$ in some cases.
\bL\label{Integ}
  If $\per{p}$ is a perversity then
$
	{\displaystyle \int} \colon \lau{\Om}{*}{\per{p}}{X} \to \Hom 
	\left( 
	\lau{RC}{\per{t} -\per{p}}{*}{X},\R \right)
$
    is  differential pairing.
\eL
\pro
Consider $\sigma \colon \Delta^i \to X$ a liftable $\per{p}$-allowable 
simplex with $\sigma(\Delta) \not\subset \Sigma_{X}$ and $\om \in \lau{\Om}{i-1}{\per{q}}{X}$. It suffices to 
prove 
\begin{equation}\label{integ}
    \int_{\sigma} d \om = \int_{\partial \sigma} \om.
\end{equation}
The boundary of $\Delta$ can be written as
$\partial \Delta = \partial_{1} \Delta + \partial_{2}\Delta$
where $\partial_{1} \Delta $ (resp. $\partial_{2} \Delta $ ) 
is composed by the faces $F$  of $\Delta$ with $\sigma 
(F) \not\subset \Sigma_{X}$ (resp. $\sigma 
(F) \subset \Sigma_{X}$). This gives the decomposition (see 
\refp{descomp}):
$$
\partial \wt\Delta = \wt{\partial_{1} \Delta}+ 
\wt{\partial_{2}\Delta }+ 
\delta \wt\Delta.
$$
We have the equalities:
$$
    \int_{\sigma} d \om \stackrel{\refp{Jos}}{=\! =}
     \int_{\wt\Delta}  \wt{\sigma}^{*}\wt{d\om} \stackrel{\refp{elem}}{=\! =}
    \int_{\wt\Delta}  d\wt\sigma^{*}\wt\om\stackrel{Stokes}{=\! = \! =\! =} \int_{\partial\wt\Delta}\wt\sigma^{*}\wt\om 
    \ \ \hbox{ and } \ \
 \int_{\partial \sigma} \om \stackrel{\refp{def},\refp{Jos}}{=\! =\! =}
\int_{\wt{\partial_{1} \Delta}}\wt\sigma^{*}\wt\om .
$$
So the equality \refp{integ} becomes ${\displaystyle 
\int_{\delta \wt{
\Delta}}\wt\sigma^{*}\wt\om 
+
\int_{\wt{\partial_{2} 
\Delta}}\wt\sigma^{*}\wt\om =0}$.
We will end the proof if we show that $\wt\sigma^{*}\wt\om=0$ on $F$, where 
the face $F$ 
\begin{itemize}
    \item[-]  is a bad face or
   \item[-]  verifies $\sigma(F) \subset \Sigma_{X}$.
   \end{itemize} 
 Put $C$ the face $\mu_{\Delta} (F)$ of $\Delta$ and $S$ the
stratum of $X$ containing $\sigma (\inte (C))$. 
Notice that the condition (All) implies
\begin{equation}
    \label{all2}
    \dim C \leq \dim F + 1 - 2 - \per{t}(S)  +
    (\per{t}(S) -\per{p}(S) )= \dim F -1  - 
    \per{p}(S). 
\end{equation}
We 
have the following commutative diagram 
$$
\begin{CD}
    \inte(F) @>\wt{\sigma}>>\mathcal{L}_{X}^{-1}(S) \\
    @V\mu_{\Delta}VV @V\mathcal{L}_{X}VV\\
    \inte(C)@>\sigma>> S\\
    \end{CD}
    $$
 It suffices to prove that the vertical degree of 
$\wt\sigma^{*}\wt\om$ relatively to $\mu_{\Delta}$ is strictly lower 
than the dimension of the fibers of $\mu_{\Delta}$, that is
$$
v_{S}\left(\wt\sigma ^{*}\wt\om \right)<
\dim F - \dim C.
$$
We distinguish two cases:
\begin{itemize}
    \item[-]
When $S$ is a regular stratum the differential form $\om$ is defined 
on $S$. We have
$$
{\wt\sigma}^{*} \wt\omega  = {\wt\sigma}^{*} 
\mathcal{L}_{X}^{*}\omega  = \mu_{\Delta}^{*}\sigma^{*} \omega , 
$$
which is a basic form relatively to $\mu_{\Delta}$. So, since $F$ is a 
bad face:
$$
v_{S}\left(\wt\sigma ^{*}\wt\om \right) \leq 0  \stackrel{\refp{bad}}{<} \dim F - 
\dim C. 
$$
\item[-] When $S$ is a singular stratum, we have
$$
v_{S}\left(\wt\sigma ^{*}\wt\om \right)
\leq 
||\om||_{S}
\leq
\per{p}(S)
\stackrel{\refp{all2}}{\leq} \dim F - \dim C -1 < \dim F - \dim C.
$$
\end{itemize}
This ends the proof.\qed

The above pairing induces the pairing
\[
     \int \colon \lau{\IH}{*}{\per{p}}{X} \TO \Hom 
     \left( 
     \lau{\IH}{\per{t} - \per{p}}{*}{X},\R \right),
 \]
 (cf. Proposition \ref{iso1}) which is not an isomorphism: for a
 cone $cL$ we have  Proposition 
 \ref{cono1} and Proposition \ref{cono5}. The problem appears when 
 negative perversities are involved. For this reason we consider 
 the relative intersection homology. 
 Since the 
 integration ${\displaystyle \int}$ vanishes on $\Sigma_{X}$ then 
 \[
	\int \colon \lau{\Om}{*}{\per{p}}{X} \TO \Hom 
	\left( 
	\lau{RC}{\per{t} - \per{p}}{*}{X,\ib{X}{\per{p}}},\R \right).
    \]
    is a well defined differential operator. We obtain the de Rham 
    duality (in the direction cohomology $\mapsto$ homology):
    \bT
    \label{Rham}
    Let $X$ be an unfoldable pseudomanifold. If $\per{p}$ is a perversity then the 
       integration induces the isomorphism
       $$
       \lau{\IH}{*}{\per{p}}{X} =
       \Hom \left( \lau{\IH}{\per{t} - \per{p}}{*}{X,\ib{X}{\per{p}}};\R\right)
       $$
    \eT
    \pro
    Following Proposition \ref{iso2} it suffices to prove that the 
    pairing
    \[
	   \int \colon \lau{\Om}{*}{\per{p}}{X} \TO \Hom 
	   \left( 
	   \lau{RC}{\per{t} - \per{p}}{*}{X,\ib{X}{\per{p}}},\R \right).
       \]
       induces an isomorphism in cohomology.
       We proceed by induction on the depth. If $\depth X = 0$ then 
       $\ib{X}{\per{p}} = \emptyset$ and we have the usual de Rham 
       theorem.
       In the general case, 
       we use the Bredon's trick (see the proof of Proposition \ref{part}) and we 
       reduce the problem to a chart 
       $X = \R^n \times cL_{S}$ with $\ib{X}{\per{p}} = \R^n \times 
       \ib{(cL_{S})}{\per{p}}$.
       Then we apply the  product formula and we reduce the problem 
       to 
    $(X ,\ib{X}{\per{p}}) = \left(cL_S,\ib{(cL_S)}{\per{p}}\right)$. We end the proof by applying 
    Corollary \ref{cono2}, Proposition 
    \ref{cono5} and the induction hypothesis.
    \qed
    
    In particular,  we have the de Rham isomorphism 
     $
	    \lau{\IH}{*}{\per{p}}{X} =
	    \lau{\IH}{\per{t} -\per{p}}{*}{X}
	    $
	    when $\per{p} \geq 
     \per{0}$.

    The intersection cohomology can be expressed in terms of the usual 
    cohomology $\coho{H}{*}{-}$ in some cases (see \cite{Bry}).
    
   \bP
   \label{XXX}
   Let $X$ be an unfoldable pseudomanifold. Then we have
   \begin{itemize}
       \item  $\lau{\IH}{*}{\per{p}}{X}= \coho{H}{*}{X - \Sigma_{X}}$ if 
       $\per{p} >\per{t}$, and
	      
       \item  $\lau{\IH}{*}{\per{q}}{X} = 
       \coho{H}{*}{X,\ib{X}{\per{q}}}$ if 
       $\per{q} \leq \per{0}$ and $X$ is normal. 
\end{itemize}
   \eP
   \pro
   From the above Theorem it suffices to prove that
       $\lau{\IH}{\per{t} -\per{p}}{*}{X}= \homo{H}{*}{X - 
       \Sigma_{X}}$ 
       and 
       \begin{equation}
	   \label{ZZZ}
	   \lau{\IH}{\per{t} -\per{q}}{*}{X,\ib{X}{\per{q}}}= 
       \homo{H}{*}{X,\ib{X}{\per{q}}}.
       \end{equation}
       The first assertion comes directly from Proposition 
       \ref{part}.  For the second one, we consider
    the differential morphism
       $$
	 A \colon \fracc{\left(\lau{AC}{\per{t} -\per{q}}{*}{X}+
	 \lau{AC}{\per{t} -\per{q}+ \per{1}}{*}{\ib{X}{\per{q}}}\right)
       \cap \partial^{-1}\left(\lau{AC}{\per{t} -\per{q}}{*-1}{X} 
       +\lau{AC}{\per{t} -\per{q}+ \per{1}}{*-1}{\ib{X}{\per{q}}}\right)}
       { \lau{AC}{\per{t} -\per{q}+ \per{1}}{*}{\ib{X}{\per{q}}}
       \cap \partial^{-1}\left(
       \lau{AC}{\per{t} -\per{q}+ \per{1}}{*-1}{\ib{X}{\per{q}}} 
       \right)}\TO 
     \frac{\homo{S}{*}{X}}{\homo{S}{*}{\ib{X}{\per{q}}}}
       $$
       defined by $A\{ \xi \} = \{ \xi \}$.
   We prove, by induction on the depth, that the morphism $A$ is a 
  quasi-isomorphism. When the 
   depth of $X$ is 0 then 
  $A$ is the identity. In the general case, 
we use the Bredon's trick (see the proof of Proposition \ref{part}) 
and we reduce the problem to a chart 
$X = \R^n \times cL_{S}$ with $\ib{X}{\per{q}} = \R^n \times 
       \ib{(cL_{S})}{\per{q}}$.
       Then we apply the  product formula and we reduce the problem to
    $\left(cL_S,\ib{(cL_S)}{\per{q}}\right)$. 
    We have three cases:
    \begin{itemize}
	   \item  $\per{q} (\vartheta)< 0$ .
		      Then
			    $\ib{(cL_{S})}{\per{q}}
			    =c\ib{(L_{S})}{\per{q}} \ne \emptyset$ and we have
		      $$
		      \lau{\IH}{\per{t} -\per{q}}{*}{cL_{S},\ib{(cL_{S})}{\per{q}}}
		      \stackrel{\ref {cono2}}{=}
		      0
		      =\homo{H}{*}{cL_{S},c\ib{(L_{S})}{\per{q}}}
		      =\homo{H}{*}{cL_{S},\ib{(cL_{S})}{\per{q}}}.
		      $$
	   
	   \item  $\per{q} (\vartheta) = 0$ and $\per{q} \ne \per{0}$ on $L_{S}$.  
		      Then
		      $\ib{(cL_{S})}{\per{q}}=c\ib{(L_{S})}{\per{q}}
		      \ne \emptyset$ and we have
		      $$
	 \lau{\IH}{\per{t} -\per{q}}{*}{cL_{S},\ib{(cL_{S})}{\per{q}}}
				\stackrel{\ref {cono2}}{=}
	\lau{\IH}{\per{t} -\per{q}}{0}{L_{S},\ib{(L_{S})}{\per{q}}}
				 \stackrel{ind}{=}
		 \homo{H}{0}{L_{S},\ib{(L_{S})}{\per{q}}}\\
				 \stackrel{norm}{=}
				 0=
			 \homo{H}{*}{cL_{S},\ib{(cL_{S})}{\per{q}}}.
				 	 $$
       
	   \item  $\per{q} =0$.    
	   Then
	   $\ib{(cL_{S})}{\per{q}}= \ib{(L_{S})}{\per{q}}  = \emptyset$ 
	   and we have
	   $$
		      \lau{\IH}{\per{t} -\per{q}}{*}{cL_{S},\ib{(cL_{S})}{\per{q}}}
					     \stackrel{\ref {cono2}}{=}
		     \lau{\IH}{\per{t} -\per{q}}{0}{L_{S},\ib{(L_{S})}{\per{q}}}
					      \stackrel{ind}{=}
	  \homo{H}{0}{L_{S},\ib{(L_{S})}{\per{q}}}\\
			    \stackrel{norm}{=}
					      \R=
		     \homo{H}{*}{cL_{S},\ib{(cL_{S})}{\per{q}}}.
						      $$
	 						    
	  \end{itemize}
	  This ends the proof. 
    \qed
   
   \prgg {\bf Remark.} Notice that we can replace the normality of $X$ 
by the connectedness of the links $\{ L_{S }\ / \ \per{q}(S) = 0\}$.

\medskip

  In the direction homology $\mapsto$ cohomology we have the following de 
  Rham Theorem
  \bC
  \label{Rhambis}
  Let $X$ be a normal unfoldable pseudomanifold. If $\per{p}$  is a  perversity then we have the isomorphism
	 $$
	 \lau{\IH}{\per{p}}{*}{X} =
	 \lau{\IH}{*}{\max ( \per{0},\per{t} -\per{p})}{X},
	 $$
  \eC
  \pro Since $\ib{X}{\max ( \per{0},\per{t} -\per{p})} = \emptyset$ 
  then 
  $
  \lau{\IH}{*}{\max ( \per{0},\per{t} -\per{p})}{X}
  $
  is isomorphic to
	  $\lau{\IH}{\per{t} -\max ( \per{0},\per{t} -\per{p})}{*}{X} 
	  =  \lau{\IH}{\min ( \per{p},\per{t} )}{*}{X}
	  $
	  (cf. Theorem \ref{Rham}). It suffices to prove that the inclusion
	  $\lau{SC}{\min ( \per{p},\per{t} )}{*}{X} \hookrightarrow
	  \lau{SC}{\per{p}}{*}{X}$
	  induces an isomorphism in cohomology.
	  We proceed by induction on the depth. When the 
	     depth of $X$ is 0 then 
	     $\lau{SC}{\min ( \per{p},\per{t} )}{*}{X} =
		       \lau{SC}{\per{p}}{*}{X} = \homo{S}{*}{X}$. In the general case, 
we use the Bredon's trick (see the proof of Proposition \ref{part}) and
we reduce the problem to 
a chart 
$X = \R^n \times cL_{S}$. We apply the  product formula     and we 
reduce the problem to 
	     to $X = cL_{S}$. Now, 
	      we have two cases
	      \begin{itemize}
		     \item  $\per{t}(\vartheta) < \per{p}(\vartheta)$. 
		   Then 
		     $
		     \lau{\IH}{\min ( \per{p},\per{t} )}{*}{cL_{S}} 
		     \stackrel{\ref{cono1}}{=}
		     \lau{\IH}{\min ( \per{p},\per{t} )}{0}{L_{S}} 
		     \stackrel{ind}{=}
		     \lau{\IH}{\per{p}}{0}{L_{S}} 
	\stackrel{\ref{cono1}, \ref{norm}}{= \! = \! = \! = \! = }
	\lau{\IH}{\per{p} }{*}{cL_{S}}.
		     $
		 
		     \item  $\per{t}(\vartheta) \geq \per{p}(\vartheta)$. 
					Then
					  $
					  \lau{\IH}{\min ( \per{p},\per{t} )}{*}{cL_{S}} \stackrel{\ref{cono1}}{=}
					  \lau{\IH}{\min ( \per{p},\per{t} )}{\leq \per{t}(\vartheta) - \per{p}(\vartheta) }{L_{S}} 
					  \stackrel{ind}{=}\lau{\IH}{\per{p}}{\leq \per{t}(\vartheta) - 
					  \per{p}(\vartheta)}{L_{S}} 
					  \stackrel{\ref{cono1}}{=}
					  \lau{\IH}{\per{p}}{*}{cL_{S}}.
					  $
		    \end{itemize}
	  This ends the proof.
  \qed
  
  \medskip

  \prgg {\bf Remark}. Notice that we can replace the normality of $X$ 
  by the connectedness of the links $\{ L_{S }\ / \ \per{p}(S) > \per{t}(S)\}$.
  In particular,  we have the de Rham isomorphism 
  $
	 \lau{\IH}{\per{p}}{*}{X} =
	 \lau{\IH}{*}{\per{t} -\per{p}}{X}
	 $
	 when $\per{p} \leq 
  \per{t}$.

    \prg {\bf Poincaré Duality}. 
    The intersection homology was introduced with the purpose of 
    extending the Poincaré Duality to singular manifolds 
    (see \cite{GM1}). The pairing is given  by the intersection 
    of cycles. For manifolds the Poincaré Duality also derives 
    from the integration of the wedge product of differential forms.
This is also the case for stratified pseudomanifolds.

    Let consider a  compact and {\em orientable} stratified pseudomanifold $X$,
    that is, the 
    manifold $X-\Sigma_{X}$ is an orientable manifold.  Let $m$ be
    the dimension of $X$.  It has been 
    proved in \cite{Bry} (see also \cite{S}) that, for a
     perversity $\per{p}$, with $\per{0} \leq \per{p} \leq \per{t}$, the pairing 
    $
  P \colon \lau{\Om}{i}{\per{p}}{X} \times 
  \lau{\Om}{m-i}{\per{t} -\per{p}}{X} \TO \R,
  $
  defined by 
  $
  P(\alpha , \beta) = {\displaystyle \int}_{X-\Sigma_{X}} \alpha \wedge \beta,
  $
  induces the isomorphism
  $
  \lau{\IH}{*}{\per{p}}{X} =\lau{\IH}{m-*}{\per{t} -\per{p}}{X}.
  $
  The same proof works for any perversity. For example, if $\per{p} < 
  \per{0}$ or $\per{p} > \per{t}$, we obtain the Lefschetz
  Duality $
  \coho{H}{*}{X,\Sigma_{X}} = \coho{H}{m-*}{X - \Sigma_{X}}
  $
(cf. Proposition \ref{XXX} and Remark 3.2.4).

\end{document}